\documentclass[reqno]{amsart}
\usepackage{float}
\usepackage{amsmath}
\usepackage{graphicx}
\usepackage{latexsym}
\usepackage{amsfonts}
\usepackage{amssymb}
\theoremstyle{plain}
\newtheorem{theorem}{Theorem}

\theoremstyle{definition}

\newtheorem*{remark*}{Remark}

\begin{document}
\title[Subcritical multitype BPRE]{Subcritical multitype branching process in random environment}

\author[Vatutin]{Vladimir Vatutin} 
\address{Novosibirsk State University and Steklov Mathematical Institute}
\email{vatutin@mi.ras.ru}

\author[Wachtel]{Vitali Wachtel} 
\address{Institut f\"ur Mathematik, Universit\"at Augsburg}
\email{vitali.wachtel@mathematik.uni-augsburg.de}

\begin{abstract}
We study the asymptotic behaviour of the survival probability of a multitype branching
process in random environment. The class of processes we consider here corresponds, in
the one-dimensional situation, to the strongly subcritical case. We also prove a
conditional limit theorem describing the distribution of the number of particles in the
process given its survival for a long time.
\end{abstract}
\maketitle
\section{Introduction and statement of results}

We study in this note asymptotic properties of multitype branching processes
in random environments (BPRE's). Informally speaking, these processes serve
as a stochastic model for the evolution of a population with $p$ different
types of particles, where each particle may have children of all types, the
number of descendants is random with distribution changing from generation to
generation in a random fashion.

In order to give a rigorous definition of multitype BPRE's and to formulate
our results we need to introduce some notation. First, $p$-dimensional
deterministic vectors with non-negative coordinates and $p\times p$
deterministic matrices with non-negative entries will be denoted by bold
lower case symbols. In particular, $\mathbf{0}:=(0,0,\ldots ,0)$ and $%
\mathbf{1}:=(1,1,\ldots ,1)$. The standard basis vectors will be denoted by $%
\mathbf{e}_{i}$, $i=1,2,\ldots ,p$. For $p$-dimensional vectors $\mathbf{x}%
=\left( x^{1},\ldots ,x^{p}\right) $ and $\mathbf{y}=\left( y^{1},\ldots
,y^{p}\right) $ we set
\begin{equation*}
(\mathbf{x},\mathbf{y}):=\sum_{i=1}^{p}x^{i}y^{i},\quad |\mathbf{x}%
|:=\sum_{i=1}^{p}|x^{i}|\quad \text{and}\quad \mathbf{x}^{\mathbf{y}%
}:=\prod_{i=1}^{p}(x^{i})^{y^{i}}.
\end{equation*}%
For every $p$-tuple $(\mu ^{1},\mu ^{2},\ldots ,\mu ^{p})$ of measures on $%
\mathbb{Z}_{+}^{p}$ we define its multidimensional generating function $%
\mathbf{f}=(f^{1},f^{2},\ldots ,f^{p})$ by the relations
\begin{equation*}
f^{i}(\mathbf{s}):=\sum_{\mathbf{z}\in \mathbb{Z}_{+}^{p}}\mathbf{s}^{%
\mathbf{z}}\mu ^{i}\{\mathbf{z}\},\quad \mathbf{s}\in \lbrack 0,1]^{p},\
i=1,2,\ldots ,p.
\end{equation*}%
Any sequence $\{\mathbf{f}_{n},n\geq 1\}$ of multidimensional generating
functions will be called a \textit{varying environment}. The corresponding
measures $\mu _{n}^{i}$ will be interpreted as the offspring law for a
particle of type $i$ in generation $n-1$. We now define a $p$-type branching
process $\mathbf{Z}_{n}=(Z_{n}^{1},Z_{n}^{2},\ldots ,Z_{n}^{p})$, $n\geq 0,$
where $Z_{n}^{i}$ is the number of type $i$ particles in the process at
moment $n$. This process has a deterministic starting point $\mathbf{Z}_{0}$
and the population sizes of the subsequent generations of the process are
specified by the following standard recursion:
\begin{equation}
Z_{n}^{i}=\sum_{j=1}^{p}\sum_{k=1}^{Z_{n-1}^{j}}X_{n,j,k}^{i},\quad
i=1,2,\ldots ,p,\ n\geq 1,  \label{def.varying}
\end{equation}%
where $\mathbf{X}_{n,j,k}=(X_{n,j,k}^{1},X_{n,j,k}^{2},\ldots
,X_{n,j,k}^{p}) $, ${k\geq 1}$ are independent random vectors distributed
according to $\mu _{n}^{j}$. Here and in what follows we denote random
objects by upper case symbols using the respective bold symbols if the
objects are vectors or matrices.

Equipping the set of all tuples $(\mu ^{1},\mu ^{2},\ldots ,\mu ^{p})$ with
the total variation distance we obtain a metric space. Therefore, we may
consider probability measures on this space. Due to the one-to-one
correspondence between $p$-tuples of measures and $p$-dimensional generating
functions, we may assume that we are given a probability measure on the
space of generating functions. With this agreement in hands we call a
sequence of independent, identically distributed random $p-$dimensional
generating functions $\{\mathbf{F}_{n}=(F_{n}^{1},F_{n}^{2},\ldots
,F_{n}^{p})\}$ a \textit{random environment}. We say that $\{\mathbf{Z}%
_{n},n\geq 0\}$ is a $p-$type BPRE, if its conditional distribution is determined by %
\eqref{def.varying} for every fixed realization of the
environmental sequence.

Branching processes in random environment with one type of particles have
been intensively investigated during the last two decades and their
properties are well understood and described in the literature. The reader
may find a modern and unified presentation of the corresponding results in a
recent book by Kersting and Vatutin~\cite{KV2017}. The multi-dimensional case
is much less studied and many basic questions, e.g. a detailed classification,
asymptotics of the survival probability and the corresponding conditional limit theorems,
are still not answered in full generality. For instance, only recently an asymptotic
representation for the survival probability in the critical multitype BPRE's
was found under relatively general conditions, see Le Page, Peigne and Pham~%
\cite{LPP2016} and Vatutin and Dyakonova~\cite{VD2016}.

The purpose of this note is to study asymptotic properties of a class of
subcritical multi-type BPRE's, which correspond to the so-called strongly
subcritical BPRE's with one type of particles.

It is a general phenomenon that the asymptotic behaviour of BPRE's is mainly
specified by the properties of certain basic characteristics of the
environment. In the case of multi-type processes the crucial role is
played by the (random) mean matrices
\begin{equation*}
\mathbf{M}_{n}=\left( M_{n}^{i,j}\right) _{i,j=1}^{p}:=\left( \frac{\partial
F_{n}^{i}}{\partial s^{j}}(\mathbf{1})\right) _{i,j=1}^{p},\quad n\geq 1.
\end{equation*}%
If $\mathbf{F}_{n}$ are independent and distributed as a generating function
$\mathbf{F}=\left( F^{1},\ldots ,F^{p}\right) $ then, obviously, $\mathbf{M}%
_{n}$ are independent probabilistic copies of the random matrix
\begin{equation*}
\mathbf{M=}\left( M^{i,j}\right) _{i,j=1}^{p}:=\left( \frac{\partial F^{i}}{%
\partial s^{j}}(\mathbf{1})\right) _{i,j=1}^{p}.
\end{equation*}
We assume that the distribution of $\mathbf{M}$ satisfies the following
assumptions:

\begin{itemize}
\item \textbf{Condition }$\mathbf{H0}.$ $\mathbb{P}\left(\left\Vert \mathbf{M%
}\right\Vert >0\right) =1$, where $\|\mathbf{M}\|$ is the operator norm of $%
\mathbf{M}$.

\item \textbf{Condition }$\mathbf{H1}$. The set $\Theta :=\left\{ \theta >0:%
\mathbb{E}\left[ \left\Vert \mathbf{M}\right\Vert ^{\theta }\right] <\infty
\right\}$ is non-empty.

\item \textbf{Condition }$\mathbf{H2}$. The support of the distribution of $%
\mathbf{M}$ acts strongly irreducibly on the semi-group of matrices with
non-negative entries, i.e. no proper finite union of subspaces of $\mathbb{R}%
^{p}$ is invariant with respect to all elements of the multiplicative
semi-group generated by the support of $\mathbf{M}$.

\item \textbf{Condition }$\mathbf{H3}$. There exists a positive number $%
\gamma >1$ such that
\begin{equation*}
1\leq \frac{\max_{i,j}M^{i,j}}{\min_{i,j}M^{i,j}}\leq \gamma .
\end{equation*}
\end{itemize}

We also need to consider the so-called Hessian matrices
\begin{equation*}
\mathbf{B}(k):=\left( \frac{\partial ^{2}F^{k}}{\partial s^{i}\partial s^{j}}%
(\mathbf{1})\right) _{i,j=1}^{p}\text{ and }\mathbf{B}_{n}(k):=\left( \frac{%
\partial ^{2}F_{n}^{k}}{\partial s^{i}\partial s^{j}}(\mathbf{1})\right)
_{i,j=1}^{p},\quad k=1,2,\ldots ,p
\end{equation*}%
and the random variables
\begin{equation*}
\mathbf{T}:=\frac{1}{\Vert M\Vert ^{2}}\sum_{k=1}^{p}\Vert \mathbf{B}%
(k)\Vert \text{ and }\mathbf{T}_{n}:=\frac{1}{\Vert M_{n}\Vert ^{2}}%
\sum_{k=1}^{p}\Vert \mathbf{B}_{n}(k)\Vert ,\quad n=1,2,\ldots .
\end{equation*}%
Thus, $\mathbf{T}_{n}$ are independent probabilistic copies of $\mathbf{T}$.
We shall impose, along with Conditions $\mathbf{H0}-\mathbf{H3}$ the following
restriction on the distribution of $\mathbf{T}:$

\begin{itemize}
\item \textbf{Condition }$\mathbf{H4}$. There exists an $\varepsilon >0$
such that
\begin{equation*}
\mathbb{E}\left[\Vert \mathbf{M}\Vert\,|\log \mathbf{T}|^{1+\varepsilon }\right] <\infty .
\end{equation*}
\end{itemize}

Using the standard subadditivity arguments, one can easily infer that for
every $\theta \in \Theta $ the limit
\begin{equation*}
\lambda \left( \theta \right) :=\lim_{n\rightarrow \infty }\left( \mathbb{E}%
\left[ \left\Vert \mathbf{M}_{n}\cdots \mathbf{M}_{1}\right\Vert
^{\theta }\right] \right) ^{1/n}<\infty
\end{equation*}%
is well defined. This function is an analogue of the moment generating
function for the associated random walk in the case of BPRE's with single
type of particles.

Set
\begin{equation*}
\Lambda (\theta ):=\log \lambda (\theta ),\quad \theta \in \Theta .
\end{equation*}

\begin{theorem}
\label{thm:nonextinct} Assume that Conditions $\mathbf{H0}-\mathbf{H4}$ are
valid, the point $\theta =1$ belongs to the interior of the set $\Theta $
and $\Lambda ^{\prime }(1)<0$. Then
\begin{itemize}
 \item[(a)] there exists a vector $\mathbf{c}=(c^1,\ldots,c^p)$ with strictly positive components such
that
\begin{equation}
\mathbb{P}\left( |\mathbf{Z}_{n}|>0\Big|\mathbf{Z}_{0}=\mathbf{e}_{i}\right)
\sim c^i\lambda ^{n}(1),\quad n\rightarrow \infty .
\label{P_asymp}
\end{equation}

\item[(b)]
for each $\mathbf{s}\in \left[ 0,1\right] ^{p}$ , $\mathbf{s}\neq \mathbf{1}$%
\begin{equation*}
\lim_{n\rightarrow \infty }\mathbb{E}\left[ \mathbf{s}^{\mathbf{Z}_{n}}\,\Big|%
|\mathbf{Z}_{n}|>0;\mathbf{Z}_{0}=\mathbf{e}_{i}\right] =\Phi _{i}\left(
\mathbf{s}\right) ,
\end{equation*}%
where $\Phi _{i}\left( \mathbf{s}\right) $ is the probability generating
function of a proper distribution on $\mathbb{Z}_{+}^{p}$.
\end{itemize}
\end{theorem}

If $p=1$ then the assumption $\Lambda ^{\prime }(1)<0$ reduces to $\mathbf{E}%
[\mathbf{M}\log \mathbf{M}]<0$. The one-type BPRE's satisfying this
condition are usually called \textit{strongly subcritical}. The asymptotic
behaviour of the survival probability for one-type strongly subcritical
processes has been studied by Guivarch and Liu~\cite{GL2001}. The
corresponding conditional limit theorem for the distribution of the number
of particles given survival has been proved by Geiger, Kersting and Vatutin
\cite{GKV2003}.

Dyakonova \cite{Dy08} has obtained \eqref{P_asymp} under the assumption that
all possible realizations of $\mathbf{M}$ have a common deterministic right
eigen-vector corresponding to the Perron root $\varrho(\mathbf{M})$ of $%
\mathbf{M}$. In this special case the condition $\Lambda^{\prime}(1)<0$
reduces to the inequality $\mathbb{E}[\varrho(\mathbf{M})\log\varrho(\mathbf{M})]<0$. The
case $\mathbb{E}[\varrho(\mathbf{M})\log\varrho(\mathbf{M})]=0$ has been
studied by Dyakonova~\cite{Dy2013} under the assumption that there exists a
common left eigen-vector corresponding to the Perron root of $\mathbf{M}$.

\section{Proof of Theorem~\protect\ref{thm:nonextinct}}

\subsection{Representation for the generating function of the process}

For every fixed realization of the environmental sequence $\mathbf{F}_{n}$ and
$k<n$ we define
\begin{eqnarray*}
\mathbf{F}_{k,n}(\mathbf{s}) &:&=\mathbf{F}_{k+1}\circ \mathbf{F}_{k+2}\circ
\cdots \circ \mathbf{F}_{n}(\mathbf{s}), \\
\mathbf{F}_{n,k}(\mathbf{s}) &:&=\mathbf{F}_{n}\circ \mathbf{F}_{n-1}\circ
\cdots \circ \mathbf{F}_{k+1}(\mathbf{s}),
\end{eqnarray*}%
and set
\begin{equation*}
\mathbf{F}_{n,n}(\mathbf{s}):=\mathbf{s}.
\end{equation*}%
It is immediate from the definition of the process $\{\mathbf{Z}_{n}, n\geq 0\}$ that
\begin{equation*}
\mathbb{E}\left[ \mathbf{s}^{\mathbf{Z}_{n}}\Big|\mathbf{Z}_{0}=\mathbf{e}%
_{i},\mathbf{F}_{1},\mathbf{F}_{2},\ldots ,\mathbf{F}_{n}\right]
=F_{0,n}^{i}(\mathbf{s}).
\end{equation*}%
Therefore,
\begin{equation}
\mathbb{E}[\mathbf{s}^{\mathbf{Z}_{n}}|\mathbf{Z}_{0}=\mathbf{e}_{i}]=%
\mathbb{E}[F_{0,n}^{i}(\mathbf{s})].  \label{start.point.1}
\end{equation}%
Setting here $\mathbf{s}=\mathbf{0}$, we get
\begin{equation}
\mathbb{P}\left( |\mathbf{Z}_{n}|>0\big|\mathbf{Z}_{0}=\mathbf{e}_{i}\right)
=1-\mathbb{E}[F_{0,n}^{i}(\mathbf{0})]=\mathbb{E}[1-F_{0,n}^{i}(\mathbf{0})].
\label{start.point.2}
\end{equation}

Let $\mathbf{f}$ be the generating function of a $p$-tuple $%
(\mu^1,\mu^2,\ldots,\mu^p)$, and let $\mathbf{m}$ be the corresponding mean
matrix, i.e.,
\begin{equation*}
\mathbf{m}=\left(\frac{\partial f^i_n}{\partial s^j}(\mathbf{1}%
)\right)_{i,j=1}^p.
\end{equation*}
For a generating function $\mathbf{f}$ and a matrix $\mathbf{a}$ define
\begin{equation*}
\psi _{\mathbf{f,a}}(\mathbf{s}):=\frac{\left\vert \mathbf{a}\right\vert} {|%
\mathbf{a}\left( \mathbf{1}-\mathbf{f}\left( \mathbf{s}\right) \right) |}-
\frac{\left\vert \mathbf{a}\right\vert }{|\mathbf{am}\left( \mathbf{1}-%
\mathbf{s}\right)|},
\end{equation*}
where, by a slight abuse of notation, $|\cdot|$ is used to denote the $L_1$%
-norm of matrices.

Let $\mathbf{a}_{i}$ be the matrix with $a^{i,i}=1$ and $a^{k,l}=0$ for all $%
(k,l)\neq (i,i)$. Then, clearly,
\begin{equation*}
1-F_{0,n}^{i}(\mathbf{s})=|\mathbf{a}_{i}(\mathbf{1}-\mathbf{F}_{0,n}(%
\mathbf{s}))|.
\end{equation*}%
Using now the definition of $\psi $, we have
\begin{align*}
\frac{1}{1-F_{0,n}^{i}(\mathbf{s})}& =\frac{|\mathbf{a}_{i}|}{|\mathbf{a}%
_{i}(\mathbf{1}-\mathbf{F}_{0,n}(\mathbf{s}))|} \\
& =\frac{|\mathbf{a}_{i}|}{|\mathbf{a}_{i}\mathbf{M}_{1}(\mathbf{1}-\mathbf{F%
}_{1,n}(\mathbf{s}))|}+\psi _{\mathbf{F}_{1},\mathbf{a}_{i}}\left( \mathbf{F}%
_{1,n}(\mathbf{s})\right)  \\
& =\frac{1}{|\mathbf{a}_{i}\mathbf{M}_{1}(\mathbf{1}-\mathbf{F}_{1,n}(%
\mathbf{s}))|}+\psi _{\mathbf{F}_{1},\mathbf{a}_{i}}\left( \mathbf{F}_{1,n}(%
\mathbf{s})\right).
\end{align*}%
Iterating this procedure, we obtain
\begin{equation*}
\frac{1}{1-F_{0,n}^{i}(\mathbf{s})}=\frac{1}{|\mathbf{a}_{i}\mathbf{R}_{n}(%
\mathbf{1}-\mathbf{s})|}+\sum_{k=1}^{n}\frac{1}{|\mathbf{a}_{i}\mathbf{R}%
_{k}|}\psi _{\mathbf{F}_{k},\mathbf{a}_{i}\mathbf{R}_{k-1}}\left( \mathbf{F}%
_{k,n}(\mathbf{s})\right) ,
\end{equation*}%
where
\begin{equation*}
\mathbf{R}_{0}:=\mathbf{Id}\quad \text{and}\quad \mathbf{R}_{k}:=\mathbf{M}%
_{1}\mathbf{M}_{2}\cdots \mathbf{M}_{k},\ k\geq 1.
\end{equation*}%
Now, recalling that $\mathbf{F}_{k}$ are i.i.d. random elements, we may
use the substitution $\mathbf{F}_{k}\leftrightarrow \mathbf{F}_{n-k+1}$. As
a result we get
\begin{align}
& \mathbb{E}[1-F_{0,n}^{i}(\mathbf{s})]=\mathbb{E}[1-F_{n,0}^{i}(\mathbf{s})]
\notag  \label{repr} \\
& \hspace{5mm}=\mathbb{E}\left[ \left( \frac{1}{|\mathbf{a}_{i}\mathbf{L}%
_{n,1}(\mathbf{1}-\mathbf{s})|}+\sum_{k=1}^{n}\frac{1}{|\mathbf{a}_{i}%
\mathbf{L}_{n,n-k+1}|}\psi _{\mathbf{F}_{n-k+1},\mathbf{a}_{i}\mathbf{L}%
_{n,n-k+2}}\left( \mathbf{F}_{n-k,0}(\mathbf{s})\right) \right) ^{-1}\right]
\notag \\
& \hspace{5mm}=\mathbb{E}\left[ \left( \frac{1}{|\mathbf{a}_{i}\mathbf{L}%
_{n,1}(\mathbf{1}-\mathbf{s})|}+\sum_{k=1}^{n}\frac{1}{|\mathbf{a}_{i}%
\mathbf{L}_{n,k}|}\psi _{\mathbf{F}_{k},\mathbf{a}_{i}\mathbf{L}%
_{n,k+1}}\left( \mathbf{F}_{k-1,0}(\mathbf{s})\right) \right) ^{-1}\right] ,
\end{align}%
where
\begin{equation*}
\mathbf{L}_{n,n+1}:=\mathbf{Id}\quad \text{and}\quad \mathbf{L}_{n,k}:=%
\mathbf{M}_{n}\mathbf{M}_{n-1}\cdots \mathbf{M}_{k},\ 1\leq k\leq n.
\end{equation*}

We shall see later, that the asymptotic, as $n\rightarrow \infty $,
behaviour of $1-F_{0,n}^{i}(\mathbf{s})$ will be determined by the summands
corresponding to fixed values $k$. To control $\psi _{\mathbf{F}_{k},\mathbf{%
a}_{i}\mathbf{R}_{k}}$ we shall use Lemma 5 in \cite{VD2016}: under the
assumption $\mathbf{H3}$
\begin{equation}
0\leq \psi _{\mathbf{F}_{k},\mathbf{a}_{i}\mathbf{L}_{n,k+1}}\left( \mathbf{F%
}_{k-1,0}(\mathbf{s})\right) \leq \gamma p^2 \mathbf{T}_{k}
\label{psi-bound}
\end{equation}
 for all $\mathbf{s}\in[0,1]^p,$ where the value of  $\psi _{\mathbf{F}_{k},\mathbf{a}_{i}\mathbf{L}_{n,k+1}}(\mathbf{s})$ at point $\mathbf{s}=\mathbf{1} $ is specified by continuity.

\subsection{Exponential change of measure.}

Define
\begin{equation*}
\mathbb{S}_+:=\{\mathbf{x}\in \mathbb{R}_+^p:|\mathbf{x}|=1\}.
\end{equation*}
Then, for every matrix $\mathbf{m}$ with positive entries the
vector
\begin{equation*}
\mathbf{m}\cdot \mathbf{x}:=\frac{\mathbf{m}\mathbf{x}}{|\mathbf{m}\mathbf{x}%
|}
\end{equation*}
is the projection of $\mathbf{m}\mathbf{x}$ on the set $\mathbb{S}_+$.

Denote by $\mathcal{C}\left( \mathbb{S}_{+}\right) $ the set of all
continuous functions on $\mathbb{S}_{+}$. For $\theta \in \Theta ,$ $g\in
\mathcal{C}\left( \mathbb{S}_{+}\right),$ and $\mathbf{x}\in\mathbb{S}_{+}$
define the transition operator
\begin{eqnarray*}
P_{\theta }g(\mathbf{x}):= \mathbb{E}\left[\left\vert \mathbf{Mx}\right\vert
^{\theta }g\left(\mathbf{M\cdot x}\right)\right].
\end{eqnarray*}

If Conditions $\mathbf{H1}-\mathbf{H3}$ hold, then, according to
Proposition 3.1 in \cite{BDGM2014}, $\lambda(\theta)$ is the spectral radius
of $P_\theta$ and there exist a unique strictly positive function $%
r_{\theta }\in \mathcal{C}\left( \mathbb{S}_{+}\right)$ and a unique
probability measure $l_\theta$ meeting the scaling
\begin{equation*}
\int_{\mathbb{S}_{+}}r_{\theta }(\mathbf{x})dl_{\theta }(\mathbf{x})=1
\end{equation*}%
and such that%
\begin{equation}
l_{\theta }P_{\theta }=\lambda \left( \theta \right) l_{\theta },\quad
P_{\theta }r_{\theta }=\lambda \left( \theta \right) r_{\theta }.
\label{invariance}
\end{equation}
Following \cite{CM2016}, we introduce the functions
\begin{equation}  \label{p-def}
p_{n}^{\theta}\left( \mathbf{x},\mathbf{m}\right) :=\frac{\left\vert \mathbf{%
mx}\right\vert ^{\theta }}{\lambda ^{n}\left( \theta \right) }\frac{%
r_{\theta }\left( \mathbf{m\cdot x}\right) }{r_{\theta }\left( \mathbf{x}%
\right) },\quad \mathbf{x}\in \mathbb{S}_{+}\text{.}
\end{equation}

It is easy to see that, for every $n\ge1$, every $\mathbf{x}\in \mathbb{S}%
_{+}$ and every matrix $\mathbf{m}$,
\begin{equation}  \label{consistency}
\mathbb{E}p_{n+1}^{\theta}\left(\mathbf{x},\mathbf{M}\mathbf{m}\right)=
p_{n}^{\theta}\left(\mathbf{x},\mathbf{m}\right)
\end{equation}
and, in particular,
\begin{equation}  \label{total_mass}
\mathbb{E}p_{n}^{\theta}\left(\mathbf{x},\mathbf{L}_{n,1}\right)=1.
\end{equation}
For each $n\ge1,$ denote by $\mathcal{F}_n$ the $\sigma$-algebra generated by
random elements $\mathbf{Z}_1,\mathbf{Z}_2,\ldots,\mathbf{Z}_n$ and $\mathbf{%
F}_1,\mathbf{F}_2,\ldots,\mathbf{F}_n$. Let $\mathbb{I}_A $ be the indicator of the event $A$. It follows from \eqref{total_mass}
that
\begin{equation*}
\mathbb{P}_n^\theta(A):= \mathbb{E}\left[p_{n}^{\theta}\left(\mathbf{x},%
\mathbf{L}_{n,1}\right)\mathbb{I}_A\right]
\end{equation*}
is a probability measure on $\mathcal{F}_n$. Furthermore, \eqref{consistency}
implies that  $\{\mathbb{P}_n^\theta,\,n\geq 1 \}$ is a sequence of consistent probability measures. It can be
extended to a probability measure $\mathbb{P}^\theta$ on our original
probability space $(\Omega,\mathcal{F})$.

We now take $\theta =1$ and apply the corresponding change of measure to the
representation \eqref{start.point.1}. Since $1-F_{n,0}^{i}(\mathbf{\mathbf{s}%
})$ is measurable with respect to $\mathcal{F}_{n}$, we have
\begin{align*}
\mathbb{E}[1-F_{n,0}^{i}(\mathbf{\mathbf{s}})]& =\lambda ^{n}(1)r_{1}(%
\mathbf{e}_{i})\mathbb{E}\left[ p_{n}^{1}(\mathbf{e}_{i},\mathbf{L}_{n,1})%
\frac{(1-F_{n,0}^{i}(\mathbf{\mathbf{s}}))}{|\mathbf{L}_{n,1}\mathbf{e}%
_{i}|r_{1}(\mathbf{L}_{n,1}\cdot \mathbf{e}_{i})}\right]  \\
& =\lambda ^{n}(1)r_{1}(\mathbf{e}_{i})\mathbb{E}^{1}\left[ \frac{%
1-F_{n,0}^{i}(\mathbf{\mathbf{s}})}{|\mathbf{L}_{n,1}\mathbf{e}_{i}|r_{1}(%
\mathbf{L}_{n,1}\cdot \mathbf{e}_{i})}\right] .
\end{align*}%
Applying now \eqref{repr}, recalling the definition of $\mathbf{a}_{i}$ and using the equality $|\mathbf{a}_{i}\mathbf{L}_{n,k}|=|\mathbf{a}_{i}\mathbf{L}_{n,k}|$
we obtain
\begin{equation*}
\mathbb{E}[1-F_{n,0}^{i}(\mathbf{\mathbf{s}})]=\lambda ^{n}(1)r_{1}(\mathbf{e%
}_{i})\mathbb{E}^{1}\left[ \frac{1}{r_{1}(\mathbf{L}_{n,1}\cdot \mathbf{e}%
_{i})}\frac{|\mathbf{e}_{i}\mathbf{L}%
_{n,1}(\mathbf{1}-\mathbf{s})|}{|\mathbf{L}_{n,1}\mathbf{e}_{i}|}\Xi _{n}(\mathbf{s})\right],
\end{equation*}%
where%
\begin{equation*}
\Xi _{n}(\mathbf{s}):=\left( 1+\sum_{k=1}^{n}\frac{|\mathbf{e}_{i}\mathbf{L}%
_{n,1}(\mathbf{1}-\mathbf{s})|}{|\mathbf{a}_{i}\mathbf{L}_{n,k}|}\psi _{%
\mathbf{F}_{k},\mathbf{e}_{i}\mathbf{L}_{n,k+1}}\left( \mathbf{F}_{k-1,0}(%
\mathbf{s})\right) \right) ^{-1}.
\end{equation*}%
Set, for brevity,
\begin{equation*}
\widetilde{\psi }_{n,k}(\mathbf{s}):=\psi _{\mathbf{F}_{k},\mathbf{a}_{i}%
\mathbf{L}_{n,k+1}}\left( \mathbf{F}_{k-1,0}(\mathbf{s})\right) .
\end{equation*}%
We fix some $N\geq 1$ and study the asymptotic behaviour of the expectation
\begin{equation*}
K_{n}(N;\mathbf{s}):=\mathbb{E}^{1}\left[ \frac{1}{r_{1}(\mathbf{L}%
_{n,1}\cdot \mathbf{e}_{i})}\frac{|\mathbf{%
e}_{i}\mathbf{L}_{n,1}(\mathbf{1}-\mathbf{s})|}{|\mathbf{L}_{n,1}\mathbf{e}_{i}|}\Xi _{n,N}(\mathbf{s})\right],
\end{equation*}%
where%
\begin{equation*}
\Xi _{n,N}(\mathbf{s}):=\left( 1+\sum_{k=1}^{N}\frac{|\mathbf{e}_{i}\mathbf{L%
}_{n,1}(\mathbf{1}-\mathbf{s})|}{|\mathbf{e}_{i}\mathbf{L}_{n,k}|}\widetilde{%
\psi }_{n,k}(\mathbf{s})\right) ^{-1}.
\end{equation*}%
For each $k\leq N$ we have the equality
\begin{equation*}
\frac{|\mathbf{e}%
_{i}\mathbf{L}_{n,k}|}{|\mathbf{e}_{i}\mathbf{L}_{n,1}(\mathbf{1}-\mathbf{s})|}=\frac{|\mathbf{e}_{i}\mathbf{L}_{n,N}\mathbf{L}%
_{N-1,k}|}{|\mathbf{e}_{i}\mathbf{L}_{n,N}\mathbf{L}%
_{N-1,1}(\mathbf{1}-\mathbf{s})|}.
\end{equation*}%
By Theorem 1 in Hennion~\cite{Hen1997}, there exist a sequence of random numbers  $\{\lambda _{n}(N)>0,\, n\geq 1\}$ and a tuple of random vectors $\{\mathbf{U%
}_{n}^{(N)}, \mathbf{V}_{n}^{(N)}, n\geq 1\}$ such that, as $n\to\infty$
\begin{equation*}
\frac{\mathbf{L}_{n,N}}{\lambda _{n}(N)}-\mathbf{V}_{n}^{(N)}\otimes \mathbf{%
U}_{n}^{(N)}\rightarrow 0\quad \text{a.s.}
\end{equation*}%
(Here and in what follows we agree to associate with vectors $\mathbf{v}=(v^1,\ldots,v^p)$ and $\mathbf{u}=(u^1,\ldots,u^p)$ the matrix $\mathbf{v}\otimes\mathbf{u}=(v_iu_j)_{i,j=1}^p.$)
Besides, the sequence of random vectors $\{(\mathbf{U}_{n}^{(N)},\mathbf{V}_{n}^{(N)}/|\mathbf{V}_{n}^{(N)}|), n\geq 1\}
$ weakly converges, as $n\to\infty,$ to a vector $(\mathbf{U}^{(N)},\mathbf{V}^{(N)})$. As a result,
the sequence of ratios
\begin{equation*}
\frac{|\mathbf{e}_{i}\mathbf{L}_{n,N}\mathbf{L}_{N-1,k}|}{|\mathbf{e}_{i}\mathbf{L}_{n,N}\mathbf{L}_{N-1,1}(\mathbf{1}-\mathbf{s}%
)|}=\frac{\left\vert \frac{%
\mathbf{e}_{i}\mathbf{L}_{n,N}}{\lambda _{n}(N)|\mathbf{V}_{n}^{(N)}|}%
\mathbf{L}_{N-1,k}\right\vert }{\left\vert
\frac{\mathbf{e}_{i}\mathbf{L}_{n,N}}{\lambda _{n}(N)|\mathbf{V}_{n}^{(N)}|}%
\mathbf{L}_{N-1,1}(\mathbf{1}-\mathbf{s})\right\vert },\quad k\leq N
\end{equation*}%
weakly converges, as $n\to\infty,$ to
\begin{equation*}
\frac{\left\vert \mathbf{e}%
_{i}(\mathbf{V}^{(N)}\otimes \mathbf{U}^{(N)})\mathbf{L}_{N-1,k}\right\vert }{\left\vert \mathbf{e}_{i}(\mathbf{V}^{(N)}\otimes \mathbf{U}^{(N)})%
\mathbf{L}_{N-1,1}(\mathbf{1}-\mathbf{s})\right\vert }%
,\quad k\leq N.
\end{equation*}%
By the same arguments, the vectors $(\widetilde{\psi }_{n,1}^{(N)}(\mathbf{s}%
),\ldots ,\widetilde{\psi}_{n,N}^{(N)}(\mathbf{s}))$ weakly converge, as $n\to\infty$
to a vector $(\widetilde{\psi }_{1}^{(N)}(\mathbf{s}),\ldots ,\widetilde{%
\psi }_{N}^{(N)}(\mathbf{s}))$, and $r_1(\mathbf{L}_{n,1}\cdot \mathbf{e}_i)$
 weakly converges, as $n\to\infty$ to a bounded random variable. Therefore, there exists
\begin{equation}
\lim_{n\rightarrow \infty }K_{n}(N;\mathbf{s})=:K^{i}(N;\mathbf{s})
\label{E1}
\end{equation}%
and the sequence $K^{i}(N;\mathbf{s})$ is decreasing in $N$ for each fixed $%
\mathbf{s=}\left( s^{1},\ldots,s^{p}\right) $.

Let $A(\mathbf{s})=\left\{ 1\leq j\leq p:s^{j}<1\right\} $. Setting $\Delta(\mathbf{s})
:=\min_{j\in A(\mathbf{s})}(1-s^{j})$ and using Lemma~2 in Furstenberg and
Kesten~\cite{FK1960} it is not difficult to check that
\begin{equation}
\frac{\Delta(\mathbf{s}) }{p^{2}\gamma ^{2}}\leq \frac{|\mathbf{e}%
_{i}\mathbf{L}_{n,1}(\mathbf{1}-\mathbf{s})|}{|\mathbf{L}_{n,1}\mathbf{e}_{i}|}\leq \frac{|\mathbf{e}_{i}%
\mathbf{L}_{n,1}|}{|\mathbf{L}_{n,1}\mathbf{e}_{i}|}\leq \gamma^2p^2,\quad n\geq 1.  \label{fk}
\end{equation}%
Combining this estimate with the fact that $r_{1}$ is bounded away from
zero, we obtain for some absolute constant $C=C(\mathbf{s})$ the estimates
\begin{align*}
0& \leq K_{n}(N;\mathbf{s})-\mathbb{E}^{1}\left[ \frac{1}{r_{1}(\mathbf{L}%
_{n,1}\cdot \mathbf{e}_{i})}\frac{|\mathbf{%
e}_{i}\mathbf{L}_{n,1}(\mathbf{1}-\mathbf{s})|}{|\mathbf{L}_{n,1}\mathbf{e}_{i}|}\,\Xi _{n}(\mathbf{s})\right]
\\
& \leq C\mathbb{E}^{1}\left[ \frac{\sum_{k=N+1}^{n}\frac{|\mathbf{e}_{i}%
\mathbf{L}_{n,1}(\mathbf{1}-\mathbf{s})|}{|\mathbf{e}_{i}\mathbf{L}_{n,k}|}%
\widetilde{\psi }_{k}(\mathbf{s})}{1+\sum_{k=N+1}^{n}\frac{|\mathbf{e}_{i}%
\mathbf{L}_{n,1}(\mathbf{1}-\mathbf{s})|}{|\mathbf{e}_{i}\mathbf{L}_{n,k}|}%
\widetilde{\psi }_{k}(\mathbf{s})}\right] .
\end{align*}%
Obviously,
\begin{equation*}
|\mathbf{e}_{i}\mathbf{L}_{n,1}(\mathbf{1}-\mathbf{s})|=|\mathbf{e}_{i}%
\mathbf{L}_{n,k}\mathbf{L}_{k-1,1}(\mathbf{1}-\mathbf{s})|\leq |\mathbf{e}%
_{i}\mathbf{L}_{n,k}|\Vert \mathbf{L}_{k-1,1}\Vert ,\quad k\leq n.
\end{equation*}%
From this bound and \eqref{psi-bound}, we get for some constant $C_1$
\begin{align*}
0& \leq K_{n}(N;\mathbf{s})-\mathbb{E}^{1}\left[ \frac{1}{r_{1}(\mathbf{L}%
_{n,1}\cdot \mathbf{e}_{i})}\frac{|\mathbf{%
e}_{i}\mathbf{L}_{n,1}(\mathbf{1}-\mathbf{s})|}{|\mathbf{L}_{n,1}\mathbf{e}_{i}|}\,\Xi _{n}(\mathbf{s})\right]
\\
& \leq C_1\mathbb{E}^{1}\left[ \frac{\sum_{k=N+1}^{\infty }\Vert \mathbf{L}%
_{k-1,1}\Vert \mathbf{T}_{k}}{1+\sum_{k=N+1}^{\infty }\Vert \mathbf{L}%
_{k-1,1}\Vert \mathbf{T}_{k}}\right] :=K(N).
\end{align*}%
Assume that we can show that the series
\begin{equation}
\Psi:=1+\sum_{k=1}^{\infty }\Vert \mathbf{L}_{k-1,1}\Vert \mathbf{T}_{k}<\infty
\quad \mathbb{P}^{1}-a.s.  \label{as-conv}
\end{equation}%
Then, clearly,
\begin{equation*}
K(N)\rightarrow 0\quad \text{as }N\rightarrow \infty .
\end{equation*}%
Combining this with \eqref{E1}, we see that there exists
\begin{equation*}
\lim_{n\rightarrow \infty }\mathbb{E}^{1}\left[ \frac{1}{r_{1}(\mathbf{L}%
_{n,1}\cdot \mathbf{e}_{i})}\frac{|\mathbf{%
e}_{i}\mathbf{L}_{n,1}(\mathbf{1}-\mathbf{s})|}{|\mathbf{L}_{n,1}\mathbf{e}_{i}|}\Xi _{n}(\mathbf{s})\right]
=\lim_{N\rightarrow \infty }K^{i}(N;\mathbf{s})=:\phi _{i}(\mathbf{s}).
\end{equation*}%
The limit at the right hand side exists, since the sequence $K^{i}(N;\mathbf{s})$ is
decreasing in~$N$. Taking into account \eqref{as-conv} and the estimate $\Xi _{n}(\mathbf{s})\geq \Psi^{-1}$  we conclude that
this limit is strictly positive for each $\mathbf{s}\in \left[ 0,1\right]
^{p}$ , $\mathbf{s}\neq \mathbf{1}$.

Having this result in hands we deduce that, as $n\rightarrow \infty $%
\begin{equation*}
\mathbb{P}\left( |\mathbf{Z}_{n}|>0\big|\mathbf{Z}_{0}=\mathbf{e}_{i}\right)
=\mathbb{E}[1-F_{0,n}^{i}(\mathbf{0})]\sim \lambda ^{n}(1)r_{1}(\mathbf{e}%
_{i})\phi _{i}(\mathbf{0})
\end{equation*}%
and%
\begin{eqnarray*}
\lim_{n\rightarrow \infty }\mathbb{E}\left[ \mathbf{s}^{\mathbf{Z}_{n}}||%
\mathbf{Z}_{n}|>0;\mathbf{Z}_{0}=\mathbf{e}_{i}\right]
&=&1-\lim_{n\rightarrow \infty }\frac{\mathbb{E}[1-F_{0,n}^{i}(\mathbf{s})]}{%
\mathbb{E}[1-F_{0,n}^{i}(\mathbf{0})]} \\
&=&1-\frac{\phi _{i}(\mathbf{s})}{\phi _{i}(\mathbf{0})}=:\Phi _{i}\left(
\mathbf{s}\right) .
\end{eqnarray*}

Therefore, to complete the proof of the theorem it remains to check the
validity of~\eqref{as-conv}.

 We first note that our assumption that $1$
belongs to the interior of $\Theta $ provides the finiteness of  $\mathbb{E}^{1}[\log \Vert
\mathbf{M}\Vert].$ Moreover, by Condition \textbf{H3}
\begin{equation*}
\min_{x\in \mathbb{S}_{+}}\frac{|x\mathbf{M}|}{|x|}\geq \gamma ^{-1}\Vert
\mathbf{M}\Vert .
\end{equation*}%
Thus, all the conditions of Theorem 2 in \cite{Hen1997} are valid and, therefore,
\begin{equation*}
\lim_{k\rightarrow \infty }\frac{\log\Vert L_{k,1}\Vert }{k}=\Lambda^{\prime} (1)<0\quad
\mathbb{P}^{1}-\text{a.s.}
\end{equation*}%
In particular, for every $\varepsilon >0$,
\begin{equation}
\Vert L_{k,1}\Vert =O\left( e^{-k^{1-\varepsilon /2}}\right) \quad \mathbb{P}%
^{1}-\text{a.s.}  \label{L_n_1}
\end{equation}%
By the Markov inequality,
\begin{equation*}
\mathbb{P}^{1}(\log \mathbf{T}_{k}>k^{(1+\delta )/(1+\varepsilon )})\leq
\frac{\mathbb{E}^{1}\left[|\log \mathbf{T}_{k}|^{1+\varepsilon }\right]}{k^{1+\delta }}.
\end{equation*}%
Using now the definition of the measure $\mathbb{P}^{1}$, we obtain
\begin{equation*}
\mathbb{P}^{1}(\log \mathbf{T}_{k}>k^{(1+\delta )/(1+\varepsilon )})\leq C_2%
\frac{\mathbb{E}[\Vert \mathbf{M}\Vert\, |\log \mathbf{T}|^{1+\varepsilon }]}{%
k^{1+\delta }}
\end{equation*}%
for some constant $C_2$. Taking into account Condition \textbf{H4} and applying the Borel-Cantelli lemma, we
conclude that
\begin{equation*}
\mathbf{T}_{k}=O\left( e^{k^{(1+\delta )/(1+\varepsilon )}}\right) \quad
\mathbb{P}^{1}-\text{a.s.}
\end{equation*}%
Combining this estimate with \eqref{L_n_1} and choosing $\delta $
sufficiently small, we obtain \eqref{as-conv}.

\vspace{12pt}

\textbf{Acknowledgement.}
This work was supported by the Russian Science Foundation under the grant 17-11-01173 and was
fulfilled in the Novosibirsk state university.

\end{document}